\documentclass{amsart}
\usepackage{amsfonts,amsmath}
\newcommand{\supp}{{\mathop{\mathrm{supp\,}}}}

\begin{document}

\title{Lifshitz asymptotics for Hamiltonians monotone in the randomness}
\author{Ivan Veseli\'c}
\address{TU Chemnitz, Fakult\"at f\"ur Mathematik, 09107 Chemnitz\\
Germany \& Emmy-Noether-Programme of the Deutsche Forschungsgemeinschaft}

\maketitle
\begin{abstract}
This is a note for the report  on the Oberwolfach Mini-Workshop: \emph{Multiscale and Variational Methods in
  Material Science and Quantum Theory of Solids}, published in \cite{CattoCVZ-07}. 
\end{abstract}

In various aspects of the spectral analysis of random Schr\"odinger operators
monotonicity with respect to the randomness plays a key role.
In particular, both the continuity properties and the low energy behaviour of the \emph{integrated density of states} (IDS)
are much better understood if such a monotonicity is present in the model than if not.

In this note we present Lifshitz-type bounds on the IDS for two classes of random potentials.
One of them is a slight generalisation of a model for which a Lifshitz bound was derived 
in a recent joint paper with Werner Kirsch \cite{KirschV}. The second one is a \emph{breather type potential}
which is a sum of characteristic functions of intervals. Although the second model is very simple, 
it seems that it cannot be treated by the methods of \cite{KirschV}.
The models and the proofs are motivated by well-established methods developed for so called 
\emph{alloy type potentials}. The basic notions of random Schr\"odinger operators and the IDS can be inferred e.g.~from 
\cite{CarmonaL-90,PasturF-92,Stollmann-01,KirschM-07,Veselic-06b}.
\medskip

\textbf{\S1 \ Random Schr\"odinger operators and the IDS. }
We consider Schr\"odinger operators on $L^2(\mathbb{R}^d)$ with a  random, $\mathbb{Z}^d$-ergodic potential.
More precisely, the random potential $W_\omega \colon \mathbb{R}^d\to\mathbb{R}$ is determined by an i.i.d.~family of
non-trivial, bounded random variables $\lambda_k\colon\Omega\to[\lambda_-,\lambda_+]=:J$ 
indexed by $k \in \mathbb{Z}^d$ and distributed according to the measure $\mu$,
and a jointly measurable \emph{single site potential} $u\colon J\times \mathbb{R}^d \to \mathbb{R}$. 
We assume that $\lambda_-\in\supp\mu$ and that $\sup_{\lambda \in J}|u(\lambda,\cdot)| \in\ell^1(L^p)$,
$p >\max(2,d/2)$. Under these assumptions the random potential 
\begin{equation}
\label{e-random-potential}
W_\omega(x)=\sum_{k\in\mathbb{Z}^d} u(\lambda_k(\omega),x-k)
\end{equation}
is relatively bounded with respect to the Laplacian with relative bound zero, uniformly in $\omega$.
Consequently, for a bounded $\mathbb{Z}^d$-periodic potential $W_{\rm per}$  the operators $H_{\rm per}:= -\Delta+W_{\rm per}$
and $H_\omega:= H_{\rm per}+W_\omega$ are selfadjoint on the domain of $\Delta$ and lower bounded uniformly in $\omega$.
Moreover, $(H_\omega)_\omega$  forms an \emph{ergodic family of operators}. Hence there exist a closed $\Sigma \subset\mathbb{R}$ and 
an $\Omega'\subset\Omega$ of full measure, such that for all $\omega \in \Omega'$ the spectrum of $H_\omega$ coincides with 
$\Sigma$. For $\Lambda_L:=[-L/2,L/2]^d, L\in \mathbb{N}$ define the distribution function 
$N(E):= L^{-d}\, \mathbb{E} \, \{{\rm Tr}[\chi_{]-\infty,E]}(H_\omega)\,\chi_{\Lambda_L}]\}$. This function is independent of $L$ and
is called IDS or \emph{spectral distribution function}. 
The support of the associated measure coincides with $\Sigma$. The IDS can be approximated 
in the sense of distribution functions by its finite volume analogs 
$N_\omega(E):= L^{-d} \sharp \{\text{eigenvalues of } H_\omega^L \le E\}$ almost surely.
Here $H_\omega^L$ denotes the restriction of $H_\omega$ to $\Lambda_L$ with Neumann boundary conditions.
For many types of random Hamiltonians the IDS is expected to be very "thin" near the spectral minimum $E_0:=\min\Sigma$.
More precisely I.~M.~Lif\v sic conjectured in \cite{Lifshitz-63,Lifshitz-64} an asymptotic behaviour of the form
$N(E) \sim c e^{-\tilde c (E-E_0)^{-d/2}}$ for $E-E_0$ small and positive, where $c, \tilde c$ denote some positive constants.
The spectrum near $E_0$ corresponds to very rare configurations of the randomness and $E_0$ is consequently 
called a \emph{fluctuation boundary}.
\medskip

\textbf{\S2 \ A class of potentials monotone in the randomness.}
Here we present a slight extension of the main result in \cite{KirschV}.
Assume that the potentials $u$ and $W_{\rm per}$ satisfy the following 
\medskip

{\sc Hypothesis A.} 
For any $\lambda \in J$ we have $\supp u(\lambda , \cdot) \subset \Lambda_1$
as well as  $u  (\lambda, x) \ge u  (\lambda_-, x)$ for all $ x \in \mathbb{R}^d$ .
There exist $\epsilon_1, \epsilon_2 > 0$ such that for all 
$\lambda \in [\lambda_-,\lambda_-+\epsilon_2]$ 
\begin{equation*}
\textstyle \int_{\mathbb{R}^d} dx \, u(\lambda, x) 
\ge 
\epsilon_1 \, (\lambda-\lambda_-)+\int_{\mathbb{R}^d} dx \, u(\lambda_-, x) 
\end{equation*}
and for all $\lambda \in [\lambda_-+\epsilon_2,\lambda_+]$ 
\begin{equation*}
\textstyle \int_{\mathbb{R}^d} dx \, u(\lambda, x) 
\ge 
\int_{\mathbb{R}^d} dx \, u(\lambda_-+\epsilon_2, x)
\end{equation*}
hold.
The function $\lambda\mapsto u(\lambda,x)$ is Lipschitz continuous at $\lambda_-$.
More precisely, for some $\kappa $, all $x\in \Lambda_1$ and all $\lambda \in [\lambda_-,\lambda_-+\epsilon_2]$
we have $u(\lambda,x)-u(\lambda_-,x) \le \kappa (\lambda-\lambda_-)$.
If $d\ge 2$, then for any $\lambda \in J$ the functions $u(\lambda,\cdot)$ and 
$W_{\rm per}$ are reflection symmetric with respect to all $d$ coordinate axes

Typical examples of potentials $u$ satisfing Hypothesis A are: an alloy type potential,
i.e.~$u(\lambda,x)=\lambda f(x)$ with $L_c^\infty(\Lambda_1)\ni f\ge 0$, and a breather type potential, 
i.e.~$u(\lambda,x)=f(x/\lambda )$ with $\supp f \subset \Lambda_{\lambda_-}, \lambda_->0, f\in C^1(\mathbb{R}^d\setminus \{0\})$ and $L^\infty(\mathbb{R}^d) \ni g(x):=  -x \cdot (\nabla f)(x)\ge 0$.
\medskip

{\sc Theorem B. (Lifshitz bound)} 
Under the Hypothesis A the IDS of the Schr\"odinger operator $H_\omega:=-\Delta+W_{\rm per} + W_\omega$
satisfies
\begin{equation}
\label{e-LifExp}\textstyle\lim_{E \searrow E_0} \frac{\log |\log N(E)|}{\log (E-E_0)} \le -\frac{d}{2}
\end{equation}
Thus for $E-E_0$ small and positive, asymptoticaly the bound 
$0< N(E) \le e^{-\tilde c (E-E_0)^{-d/2}}$ holds. The proof is essentially the same as in \cite{KirschV}.
\medskip

\textbf{\S3 \ Breather potentials with characteristic functions of intervals.}
We consider a very explicite class of random potentials on $\mathbb{R}$. 
Let $(\lambda_k)_{k\in \mathbb{Z}}$ be as before with $\lambda_-=0$, $\lambda_+=1$.
The breather type potential
\begin{equation}
\label{e-characteristic-breather}\textstyle
W_\omega(x)= \sum_{k\in\mathbb{Z}} u(\lambda_k(\omega),x-k), \quad
\text{ where } u(\lambda,x)= \chi_{]0,\lambda]}(x)
\end{equation}
does not satisfy the Lipschitz condition in Hypothesis A. Nevertheless we have 
\medskip

{\sc Theorem C.} 
The IDS of the Schr\"odinger operator $H_\omega:=-\Delta+W_\omega$, 
where $W_\omega$ is as in \eqref{e-characteristic-breather},
satisfies the Lifshitz bound \eqref{e-LifExp}. Note that $E_0=0$ for this model.
\medskip

It seems that the reason why the method of \cite{KirschV} is not applicable to the
potential \eqref{e-characteristic-breather} is the use of Temple's inequality \cite{Temple-28}.
For Temple's inequality to yield an efficient estimate, the second moment $\langle H_\omega^L\psi,H_\omega^L\psi\rangle$ in an well chosen state
$\psi$ has to be much smaller than the first moment $\langle \psi,H_\omega^L\psi\rangle$. 
For the current application the best choice of $\psi$ seems to be the periodic, positive ground state of 
$H_{\rm per}$. However for such $\psi$ and for the potential \eqref{e-characteristic-breather}, the first 
and second moment coincide! It turns out that Thirring's inequality \cite[3.5.32]{Thirring-94} is better adapted to 
the model under consideration. It was used before in \cite{KirschM-83a} in a similar context.
\smallskip

{\sc Sketch of proof:} 
As before the superscript ${}^L$ denotes the Neumann b.~c.~restriction to $[-L/2,L/2]$.
Since $N(E) \le L^{-1} {\rm Tr}[\chi_{]-\infty,E]}(-\Delta^L)] \mathbb{P}\{\omega\mid E_1(H_\omega^L) \le E\} $
for any $L\in \mathbb{N}$, it is sufficient to derive an exponential bound on the probability that the first eigenvalue $E_1$ of 
$H_\omega^L$ does not exceed $E$. 

We set $I_L:= \Lambda_L \cap \mathbb{Z}, H_0:= -\Delta -\alpha/4L^2$, $\psi =L^{-1/2}\chi_{\Lambda_L}$ and 
$V_\omega(x)=\alpha/4L^2+W_\omega(x)$. Then $E_1(H_0^L)=-\alpha/4L^2$ and $E_2(H_0^L)\ge 3\alpha/4L^2$, cf.~\cite{KirschS-87}. 
Since $V_\omega$ does not vanish, $V_\omega^{-1}$ is well-defined and we calculate 
$L(\int_{\Lambda_L} V_\omega(x)^{-1} dx)^{-1} =\frac{\alpha}{4L^2} \frac{4L^2+\alpha}{4L^2-4L^2S_L+\alpha}$.
We use the notation $S_L:= L^{-1} \sum_{k \in I_L} \lambda_k$ for averages, $\tilde \lambda_k:=\min(\lambda_k, 1/2)$ for cut-off random variables and similarly $\tilde V_\omega$, $\tilde S_L$ 
for the potential and the averages. Then $E_1(H_0^L) + \langle \psi, \tilde V_\omega^{-1} \psi\rangle^{-1} \le \alpha/4L^2 < E_2(H_0^L)$,
thus Thirring's  inequality is applicable and yields
\[
\textstyle 
E_1(H_\omega^L) \ge E_1(\tilde H_\omega^L) \ge E_1(H_0^L) + \langle \psi, \tilde V_\omega^{-1} \psi\rangle^{-1}
\ge \frac{\alpha \tilde S_L}{5 L^2}
\]
as soon as $L^2 \ge \alpha$. For a given $E>0$ chose $L:= \lfloor\beta E^{-1/2}\rfloor$, then 
$\mathbb{P}\{E_1(H_\omega^L) \le E\} \le \mathbb{P}\{\alpha\tilde S_L/5L^2 \le E\} \le \mathbb{P}\{\alpha\tilde S_L/5 \le \beta^2\}$.
Since $0<\, \mathbb{E} \,\{\tilde S_L\} =\, \mathbb{E} \,\{\tilde \lambda_k\} \le 1/2$ it is possible to choose $0<\beta \le \sqrt{\alpha\, \mathbb{E} \,\{\tilde \lambda_k\}/10}$.
With this choice we have 
$\mathbb{P}\{\tilde S_L \le 5\beta^2/\alpha\} \le \mathbb{P}\{\tilde S_L \le \, \mathbb{E} \,\{\tilde S_L\}/2\}$. A large deviation estimate bounds this 
probability by $ e^{-c L^d} = e^{-\tilde c E^{-d/2}}$ for some positive constants $c, \tilde c$. This completes the proof.

The higher dimensional analog of this model is currently under study.

\end{document}